\documentclass[12pt]{amsart}
\usepackage{amsmath,amssymb,amsthm}
\usepackage[utf8]{inputenc}
\usepackage[english]{babel}
\usepackage{enumitem}
\usepackage{todonotes,xcolor}
\usepackage{graphicx}
\usepackage[colorlinks=true,linkcolor=blue,citecolor=blue,urlcolor=blue]{hyperref}
\usepackage[alphabetic]{amsrefs}

\newtheorem{lemma}{Lemma}
\newtheorem{theorem}{Theorem}
\newtheorem{theoA}{Theorem}

\newtheorem{prop}{Proposition}
\newtheorem{cor}{Corollary}
\theoremstyle{remark}
\newtheorem{remark}{Remark}

\title{A note on generators of coordinate rings}
\author[Jari Desmet]{\texorpdfstring{Jari Desmet$^*$}{Jari Desmet}}
\author[Ayan Dey]{\texorpdfstring{Ayan Dey$^\dagger$}{Ayan Dey}}
\thanks{$^*$University of Ghent, Belgium. \emph{E-mail}: \texttt{Jari.Desmet@UGent.be}}
\thanks{$^\dagger$Indian Statistical Institute, New Delhi, India. \emph{E-mail}: \texttt{studentayandey@gmail.com}}

\begin{document}
\maketitle
 
\begin{abstract}
    Given a connected simply connected simple complex Lie group $G$, it is known that the coordinate ring $\mathcal{O}(G)$ over $\mathbb{C}$ coincides with the Berenstein--Fomin--Zelevinsky cluster algebra when $G$ is not of type $F_4$. In this short note, we settle the final open case, $F_4$, by showing that $\mathcal{O}(F_4)$ is generated by generalized minors over $\mathbb{C}$.
\end{abstract}

\section{Introduction}

Motivated by the theory of canonical bases, Fomin and Zelevinsky in \cite{FZ02} introduced cluster algebras to provide a systematic method for studying canonical bases of (quantized) coordinate rings. One of the major goals of this program is to show that the coordinate rings of several varieties related to a simply connected simple complex Lie group $G$ carry natural cluster algebra structures. As of now, this has been verified for a wide range of examples. For the coordinate ring $\mathcal{O}(G)$ itself, when $G$ is of type $A_n$ ($n\geq 1$), it is shown in \cite{FWZ21} that $\mathcal{O}(G)$ is a cluster algebra. The $n=1$ case was already proved in \cite{FZ02} as a concrete example of a cluster algebra. Beyond type $A_n$ ($n\geq 1$), this conjecture was recently proved by Oya in \cite{Oya} for every $G$ except type $F_4$. More specifically, Oya proved the chain of equalities $\mathcal{O}(G)=\mathcal{A}=\mathcal{U}$ when $G$ is not of type $F_4$. Here, $\mathcal{A}$ is the partially compactified cluster algebra associated with the Berenstein--Fomin--Zelevinsky seed for the open double Bruhat cell in $G$, and $\mathcal{U}$ is the corresponding upper cluster algebra. Almost in parallel, the equality $\mathcal{O}(G)=\mathcal{U}$ was established for every $G$ (including $F_4$) by Qin and Yakimov in \cite{QY}. As noted by Oya, using \cite[Theorem 1.4]{GLS13}, the proof of the equalities $\mathcal{O}(G)=\mathcal{A}=\mathcal{U}$ reduces to proving the inclusion $\mathcal{O}(G)\subseteq \mathcal{A}$, which can be obtained by proving that the coordinate ring $\mathcal{O}(G)$ is generated by generalized minors over $\mathbb{C}$. This property is almost straightforward when $G$ is of types $A_n$--$D_n$, $E_6$, and $E_7$, because $G$ admits a tensor generating family of minuscule representations. For $G$ of type $G_2$, $F_4$, or $E_8$, there are no non-trivial minuscule representations. Instead, $\mathcal{O}(G)$ is generated by matrix coefficients of the quasi-minuscule module $V_s:=V(\theta_s)$, where $\theta_s$ is the highest short root of $\mathrm{Lie}(G)$. In each of these three types, there is a unique (up to a non-zero scalar) non-trivial $\mathfrak{g}$-module morphism 
\[
m: V_s \otimes V_s \to V_s.
\]
Furthermore, for $G_2$ and $E_8$, we have $m((V_s)_0\otimes (V_s)_0)=\{0\}$ (see \cite[Theorem 3.9]{Oya}). This vanishing property is the key observation that is subsequently exploited to prove the generation via generalized minors for types $G_2$ and $E_8$. However, this technique fails for $F_4$ because $m((V_s)_0\otimes (V_s)_0)\neq \{0\}$ as pointed out in \cite[Remark 3.10]{Oya}. Therefore, the goal of the present article is to settle the final $F_4$-case; in fact, our proof is type independent, and we obtain the following generation theorem in complete generality, which may be of independent interest as well.

\begin{theoA}\label{th:main}
    Given any connected simply connected simple complex Lie group $G$, $\mathcal{O}(G)$ is generated by generalized minors over $\mathbb{C}$.
\end{theoA}

Using Theorem \ref{th:main} and the same proof as in \cite[Theorem 3.1]{Oya}, we then immediately obtain the required chain of equalities: $\mathcal{O}(G)=\mathcal{A}=\mathcal{U}$ without any assumption on $G$. This leads to new applications for studying cluster structures on $\mathcal{O}(X)$ for some variety $X$ (over $\mathbb{C}$) with a morphism $X\to G$, strengthening \cite[Theorem 1]{IOS23} as mentioned in \cite{Oya}, and also for studying cluster structures on the $(G, c)$-bands ($c$ is a Coxeter element) recently introduced by Francone and Leclerc in \cite{FL25}.

\subsection*{Strategy of the proof}\label{st}
Surprisingly, the quantum analog of the above theorem is now fully established, including type $F_4$. Indeed, when $G$ is not of type $F_4$, Oya, Qin and Yakimov in \cite{OQY} showed that the quantized coordinate ring $\mathcal{O}_q(G)$ is generated by generalized quantum minors. Their approach, specifically in the type $E_8$ case, was then extended by the second author in \cite[Theorem 5.3]{Dey} uniformly across the types $G_2$, $F_4$, and $E_8$ by studying the product crystal $B(\theta_s)\otimes B(\theta_s)$ of the tensor product module $V_s\otimes V_s$. Although the quantum approach cannot be directly specialized at $q=1$, the motivation for our present proof comes from it. For the quantized coordinate ring $\mathcal{O}_q(F_4)$, the analogous obstruction, $m((V_s)_0\otimes (V_s)_0)\neq \{0\}$, is bypassed using the braiding isomorphism $\hat{R}$ on $V_s \otimes V_s$. The quantum analog of $m$ satisfies the fundamental relation 
\[
m \circ \hat{R} = d \cdot m,\quad d\neq 1.
\]
Twisting $m$ by the braiding $m(\mathrm{id} - \hat{R}) = (1 - d)m$ then cancels out the $(V_s)^*_0 \otimes (V_s)^*_0$ component on the dual side (see \cite[Lemma 5.2]{Dey}). However, as $q \to 1$, $\hat{R}$ specializes to the flip map $\tau: v \otimes w \mapsto w \otimes v$. Because $V_s$ appears in the symmetric subspace of $V_s \otimes V_s$, $d \to 1$ as $q \to 1$. Therefore, the relation degenerates to $m \circ \tau = m$, and we get $m(v \otimes w - w \otimes v) = 0$.

This failure then naturally leads us to look for a simple module $V(\lambda)$ inside the anti-symmetric subspace of $V_s \otimes V_s$ whose projection $D: V_s \otimes V_s \to V(\lambda)$ satisfies $D \circ \hat{R} = d \cdot D$ with $d \neq 1$ such that $d$ approaches $-1$ as $q \to 1$. The adjoint module $V_l$ satisfies all such properties for type $F_4$ and occurs in $V_s \otimes V_s$ with multiplicity $1$. We shall use the properties of the map $D$ at $q=1$ such as $D((V_s)_{0}\otimes(V_s)_{0})=\{0\}$ to prove the desired generation theorem. Once this key ingredient is identified, we shall neither use quantum groups nor cluster algebra techniques to prove Theorem \ref{th:main}; therefore, the rest of the article relies entirely on pure Lie-theoretic tools.

\section{Preliminaries and Background}

\noindent We follow the standard notation of Lie theory \cite{Hum}. Let $\mathfrak{g}$ be a finite-dimensional complex simple Lie algebra with a fixed Cartan subalgebra $\mathfrak{h}$, simple roots $\Delta:= \{\alpha_i\}_{i \in I}$, root system $\Phi \subseteq \mathfrak{h}^*$, coroot system $\Phi^\vee=\{\alpha^\vee\mid \alpha\in \Phi\}\subseteq \mathfrak{h}$, weight lattice $P$, and dominant weights $P_+$. For each $\alpha\in \Phi$, we define the reflection $\sigma_\alpha\in \mathrm{GL}(\mathfrak{h})$,
\[
\sigma_\alpha(h):=h-\alpha(h)\alpha^\vee
\]
for any $h\in \mathfrak{h}$, and let $W = \langle \sigma_\alpha \mid \alpha \in \Phi \rangle$ be the Weyl group. Equivalently on $\mathfrak{h}^*$, the action of $W$ is defined by $(w.\lambda)(h):=\lambda(w^{-1}(h))$ for $\lambda\in \mathfrak{h}^*$, $w\in W$, and $h\in \mathfrak{h}$. We follow the ordering of the Dynkin diagrams as given in \cite[Theorem 11.4]{Hum}.
\vskip 0.1cm
Recall that the roots of $\mathfrak{g}$ can have at most two different lengths. When all the roots of $\mathfrak{g}$ are of the same length, we call such a $\mathfrak{g}$ simply laced and these are the types $A_n$ ($n\geq 1$), $D_n$ ($n\geq 4$), $E_6$, $E_7$, and $E_8$. Furthermore, it is customary to refer to all the roots as short in this case. Otherwise, if $\mathfrak{g}$ is non-simply laced, that is, of type $B_n$ ($n\geq 2$), $C_n$ ($n\geq 3$), $F_4$, or $G_2$, it has roots of two different lengths, classified as short roots and long roots. We denote by $\Phi_s$ and $\Phi_l$ the sets of short and long roots in $\Phi$, respectively. Therefore, there exists a unique $W$-invariant positive definite symmetric bilinear form $(\;,\;)$ on $\mathbb{R}.\Phi$ normalized so that $(\alpha, \alpha) = 2$ for all short roots $\alpha$.

The universal enveloping algebra of $\mathfrak{g}$ is denoted by $U(\mathfrak{g})$. Its Hopf structure is defined using the maps:
\[
\Delta(x) = x \otimes 1 + 1 \otimes x, \quad \varepsilon(x) = 0, \quad S(x) = -x  
\]
for $x \in \mathfrak{g}$.

\subsection*{Representations of $U(\mathfrak{g})$}
Every finite-dimensional $U(\mathfrak{g})$-module $M$ decomposes into weight spaces $M = \bigoplus_{\nu \in P} M_\nu$, where $M_\nu := \{u \in M \mid h \cdot u = \nu(h)u \text{ for all } h \in \mathfrak{h}\}$, and we set $M_{\neq 0} := \bigoplus_{\nu \neq 0} M_\nu$. The dual $M^*$ is a $U(\mathfrak{g})$-module via $\langle x \cdot f, u \rangle := \langle f, S(x)u \rangle$. Simple finite-dimensional $U(\mathfrak{g})$-modules are the highest weight modules $V(\lambda)$ ($\lambda \in P_+$), whose extremal weight spaces $V(\lambda)_{w\lambda}$ ($w \in W$) are one-dimensional.

Recall that a non-trivial simple module $V(\varpi)$ is \emph{minuscule} if its weights form the single orbit $W \cdot \varpi$, and \emph{quasi-minuscule} (or \emph{little adjoint}) if its weights are $W \cdot \varpi \cup \{0\}$. The \emph{adjoint module} $\mathfrak{g} \cong V_l$ has the highest weight equal to the highest root $\theta$ of $\Phi$. We also have $(V_l)_0 = \mathfrak{h}$ and one-dimensional root spaces $(V_l)_\alpha = \mathbb{C}.e_\alpha$ ($\alpha \in \Phi$).

The unique non-minuscule \emph{quasi-minuscule} module $V_s$ has the highest weight $\theta_s\in \Phi_s$, the highest short root. When $\mathfrak{g}$ is simply laced, we have $\theta=\theta_s$, and therefore $V_l = V_s$. Also, note that the modules $V_s$ and $V_l$ are self-dual.

Let $G$ be the connected simply connected complex Lie group with $\mathrm{Lie}(G) \cong \mathfrak{g}$. Then any finite-dimensional $\mathfrak{g}$-module $M$ uniquely induces a $G$-action on $M$, and any $\mathfrak{g}$-equivariant morphism between two finite-dimensional modules is also $G$-equivariant, and vice versa. Let $T \subseteq G$ be the maximal torus with $\mathrm{Lie}(T) = \mathfrak{h}$, and let $N(T)$ be its normalizer in $G$. Under the isomorphism $W \cong N(T)/T$, the Weyl group $W$ permutes the weight spaces of $M$. That is, for any $w = gT \in W$ with $g \in N(T)$ and any $v_\mu \in M_\mu$, we have $v_{w\mu}:= g \cdot v_\mu \in M_{w\mu}$. This is well-defined up to a non-zero scalar, depending on the choice of representative $g$. In particular, if $M$ has a non-trivial weight space $M_0$ of weight 0, this defines a genuine $W$-action on $M_0$ via $w \cdot v_0:= g \cdot v_0$ as $T$ acts trivially on $M_0$. Note that if $\mathfrak{g}$ is simply laced, we have $V_s=V_l$ and $(V_s)_0=(V_l)_0=\mathfrak{h}$. So, $W$ acts irreducibly on $(V_s)_0$. For non-simply laced $\mathfrak{g}$ we have the following result. 

\begin{lemma}\label{lem:weyl-zero}
	Assume that $\mathfrak{g}$ is non-simply laced. The kernel of the action of $W$ on $(V_s)_0$ contains $\sigma_\beta$ for all long roots $\beta$, and $W_s = \langle \sigma_\alpha \mid \alpha \in \Delta \text{ short}\rangle$ acts irreducibly on $(V_s)_0$. 
\end{lemma}
\begin{proof}
	The first statement follows from \cite[Theorem~3.1]{quasimin}, while the second follows from  \cite[Lemma~3.4]{quasimin}.
\end{proof}

\subsection*{Coordinate rings}
For a finite-dimensional $U(\mathfrak{g})$-module $M$, $v \in M$, and $f \in M^*$, we write $C^M(f, v) \in U(\mathfrak{g})^*$ for the matrix coefficient $X \mapsto \langle f, X \cdot v \rangle$. The coordinate ring $\mathcal{O}(G) \subseteq U(\mathfrak{g})^*$ is spanned by the matrix coefficients of all finite-dimensional $U(\mathfrak{g})$-modules.

For $\lambda \in P_+$ and $w_1, w_2 \in W$, fix two non-zero vectors $v_{w_i\lambda} \in V(\lambda)_{w_i\lambda}$ and let $f_{-w_i\lambda} \in (V(\lambda)^*)_{-w_i\lambda}$ ($i=1,2$) be their dual functionals. The associated \emph{generalized minor} is defined to be the matrix coefficient
\[
\Delta_{w_1\lambda, w_2\lambda} := C^{V(\lambda)}(f_{-w_1\lambda}, v_{w_2\lambda}),
\]
and it is unique up to a non-zero scalar. Observe that any $C^{V_s}(g_\gamma, v_\mu)$ ($\gamma, \mu \in \Phi_s$) is a generalized minor up to a scalar.

\section{Main results}

We start with the following folklore lemma. 

\begin{lemma}\label{lem:short-roots}
	Assume $\mathfrak{g}$ is of rank at least $2$ and not of type $B_n$ ($n\geq 2$). Given any non-zero $h \in \mathfrak{h}$, there are at least two short roots $\alpha, \beta$ such that $\alpha \neq \pm \beta$ and both $\alpha(h) \neq 0$ and $\beta(h) \neq 0$.
\end{lemma}
\begin{proof}
	As short roots span $\mathfrak{h}^*$, there is at least one short root $\alpha \in \Phi_s$ such that $\alpha(h) \neq 0$. Assume, to the contrary, that for any other short root $\gamma \in \Phi_s$ with $\gamma \neq \pm \alpha$, we have $\gamma(h) = 0$. Then $\sigma_\alpha(\gamma)$ is again a short root other than $\pm \alpha$, so $\sigma_\alpha(\gamma)(h) = 0$. Since
	\[
		\sigma_\alpha(\gamma)(h) = \gamma(h) - \frac{2(\gamma,\alpha)}{(\alpha,\alpha)}\alpha(h)
	\]
	and $\alpha(h) \neq 0$, this implies $(\gamma,\alpha) = 0$ for all $\gamma \in \Phi_s$ with $\gamma \neq \pm \alpha$. As all the short roots form a single $W$-orbit of $\Phi$ and the form $(\;,\;)$ is $W$-invariant, this implies that when $\mathfrak{g}$ is at least of rank $2$ and not of type $B_n$ ($n\geq 2$) or $G_2$, for any two distinct short simple roots $\alpha_i$ and $\alpha_j$ we get $(\alpha_i, \alpha_j) = 0$. But the Dynkin diagrams (see \cite[Theorem 11.4]{Hum}) of these types show otherwise, as we can always find two distinct short simple roots connected by an edge. For $\mathfrak{g} = G_2$, take the short simple root $\alpha_1$ and the highest short root $2\alpha_1+\alpha_2$ (see Table 2 of \cite[Chapter 12]{Hum}); then we have $(\alpha_1, 2\alpha_1 + \alpha_2) \neq 0$. Therefore, when $\mathfrak{g}$ is at least of rank $2$ and not of type $B_n$ ($n\geq 2$), there must be at least two short roots $\alpha, \beta \in \Phi_s$ such that $\alpha \neq \pm \beta$ and $\alpha(h) \neq 0$, $\beta(h) \neq 0$.
\end{proof}

Given any non-trivial finite-dimensional $\mathfrak{g}$-module $M$, $V_l \;(\cong \mathfrak{g})$ is always embedded inside $\mathrm{End}(M)\cong M\otimes M^*$ as a $\mathfrak{g}$-module. So, using the self-duality of $V_s (=V(\theta_s))$, we always have  $V_l\subseteq V_s\otimes V_s$. 

More generally, $V_s\otimes V_s$ decomposes into the direct sum of all $V(\theta_s+\mu)$ with $\mu \in W\theta_s$ and $\theta_s+\mu\in P_+$, and $r$ copies of $V_s$ where $r$ is the number of short simple roots $\alpha\in \Delta_s$ with $(\theta_s, \alpha)>0$ (see \cite[Lemma 5A.9b]{Jan}). From this, it now follows that when $\mathfrak{g}$ is non-simply laced, there is a unique copy of $V_l\cong V(\theta)$ inside $V_s\otimes V_s$ because $\theta-\theta_s\in \Phi_s=W\theta_s$. Furthermore, if $\mathfrak{g}$ is non-simply laced and not of type $B_n$ ($n\geq 2$), there is a unique copy of $V_s$ inside $V_s\otimes V_s$. For $\mathfrak{g}=B_n$ ($n\geq 2$), there is no copy of $V_s$ inside $V_s\otimes V_s$ as $(\theta_s, \alpha_n)=0$ for the unique short simple root $\alpha_n\in \Delta_s$.

We use the following non-zero $\mathfrak{g}$-equivariant maps:
\begin{align*}
    m: V_s \otimes V_s &\to V_s,\\
    D: V_s\otimes V_s &\to V_l, \\
	\mathrm{ad}: V_l \otimes V_s &\to V_s. 
\end{align*}
Here, the map $\mathrm{ad}$ is defined using the $V_l\cong \mathfrak{g}$-action on the quasi-minuscule module $V_s$. If $\mathfrak{g}$ is simply laced, we have $V_l=V_s$, and therefore, we shall assume that the maps $m$, $D$, and $\mathrm{ad}$ coincide with the Lie bracket $[\;,\;]$ of $\mathfrak{g}$. From the discussion in the previous paragraph, it now follows that if $\mathfrak{g}$ is non-simply laced and not of type $B_n$ ($n\geq 2$), both maps $m$ and $D$ exist and are unique up to a non-zero scalar. Let us record some important properties of these maps.

\begin{prop}\label{prop:ad}
Let $\mathfrak{g}$ be a simple complex Lie algebra with $V_s$ and $V_l$ being its quasi-minuscule and adjoint modules, respectively.
    \begin{enumerate}
	\item For $v_\mu\in V_s$ with $\mu\in \Phi_s$, there is an element $h\in (V_l)_0$ such that $\mathrm{ad}(h\otimes v_\mu)=v_\mu$.
        \vskip 0.1cm
	\item $\mathrm{ad}((V_l)_0\otimes (V_s)_0) = \{0\}$.
    \end{enumerate}
\end{prop}

\begin{proof}
    To prove (1), recall that $\mathfrak{g}\cong V_l$ as a $\mathfrak{g}$-module. So, after normalizing if needed, there is an element $h\in\mathfrak{h}\;(= (V_l)_0)$ such that $\mathrm{ad}(h\otimes v_\mu)=v_\mu$. Similarly, (2) follows from the definition of weight-zero vectors of $(V_s)_0$.
\end{proof}

\begin{prop}\label{prop:maps}
Assume that $\mathfrak{g}$ is at least of rank 2 and not of type $B_n$ ($n\geq 2$). Then we have
\begin{enumerate}
	\item
    $\mathrm{Span}_{\mathbb{C}}\{m(v_\mu\otimes v_{-\mu})\mid \mu\in \Phi_s\}=(V_s)_0$,
    \vskip 0.1cm
	\item $D((V_s)_{0}\otimes(V_s)_{0})=\{0\}$ and $D((V_s)_{\neq 0}\otimes(V_s)_{\neq 0}) = V_l$.
\end{enumerate}
\end{prop}

When $\mathfrak{g}$ is simply laced, some of these arguments are trivial because of $[\mathfrak{h},\mathfrak{h}]=0$. To motivate the proof for the non-simply laced case, we first focus on the case when $\mathfrak{g}=\mathfrak{f}_4$. More specifically, the Lie algebra $\mathfrak{f}_4$ can be realized as the algebra of derivations $\mathrm{Der}(\mathfrak{J})\cong \mathfrak{f}_4$ of the $27$-dimensional exceptional Jordan algebra or Albert algebra $\mathfrak{J}$ over complex octonions $\mathbb{O}_{\mathbb{C}}$. Here, $\mathfrak{J}$ consists of $3\times 3$ Hermitian matrices with entries from $\mathbb{O}_{\mathbb{C}}$. Under the module action of $\mathrm{Der}(\mathfrak{J})$ on $\mathfrak{J}$, $\mathfrak{J}$ decomposes as $\mathfrak{J} \simeq V_s \oplus \mathbb{C}1$, where $V_s$ is the unique 26-dimensional quasi-minuscule module of $\mathfrak{f}_4$ consisting of trace-free matrices of $\mathfrak{J}$, and $1$ stands for the identity matrix. Under this realization, $(V_s)_0$ corresponds to the subspace of trace-free diagonal matrices in $\mathfrak{J}$ (see \cite[Section~2]{Yok} and \cite[Remark~3.10]{Oya}). For any $X, Y \in V_s \subseteq \mathfrak{J}$,
\[
    D(X \otimes Y) := [L_X, L_Y],
\]
where $L_X$ and $L_Y$ are the left multiplication operators on $\mathfrak{J}$ via the Jordan product\footnote[1]{For any two $X, Y\in\mathfrak{J}$, the Jordan product of $X$ and $Y$ is $\frac{1}{2}(XY+YX)$.}. Since $[L_X, L_Y]$ is a derivation of $\mathfrak{J}$ (see \cite[Proposition 2.4.1]{Yok}), it lies in $\mathrm{Der}(\mathfrak{J}) \cong \mathfrak{f}_4 \cong V_l$. Using the fact that $(V_s)_0$ consists of trace-free diagonal matrices in $\mathfrak{J}$, it now follows that $D((V_s)_0 \otimes (V_s)_0) = \{0\}$.  

\begin{proof}
To prove (1), we need two observations. First, the algebra $(V_s,m)$ (with the algebra structure $v.w:= m(v\otimes w)$) is simple by \cite[Theorem 5(i)]{popov}. This implies that there exist weight vectors $v_{\pm \mu} \in V_s$ of weight $\pm \mu\in \Phi_s$ such that $m(v_\mu\otimes v_{-\mu})\neq 0$. Indeed, the subspace $(V_s)_{\neq0}$ would be a proper ideal otherwise. Now, $W$ permutes the non-zero weight spaces of $V_s$, and acts irreducibly on $(V_s)_0$ by Lemma~\ref{lem:weyl-zero} (and as noted earlier in the simply laced case). Hence, $\text{Span}_{\mathbb{C}}\{m(v_{w\mu}\otimes v_{-w\mu}) \mid w\in W\} =(V_s)_0$, and this proves (1).

Finally, to prove (2) first assume that $\mathfrak{g}$ is simply laced, so we have $D=[\;,\;]$, and hence the proof follows from the fact that $[\mathfrak{h},\mathfrak{h}]=0$. 

Next, assume $\mathfrak{g}$ is non-simply laced and recall that for any long root $\beta\in \Phi_l$, the reflection $\sigma_\beta$ acts trivially on $(V_s)_0$ by Lemma~\ref{lem:weyl-zero}. Hence, $D((V_s)_0\otimes (V_s)_0)$ is a subspace of the Cartan subalgebra on which  $\sigma_\beta$ acts trivially for every long root $\beta$. Since $\Phi_l$ spans $\mathfrak{h}^*$, we have $D((V_s)_0\otimes (V_s)_0)=0$, and therefore, by comparing weights on both sides, it follows immediately that $(V_l)_0\subseteq D((V_s)_{\neq 0}\otimes\allowbreak (V_s)_{\neq 0})$ and $(V_l)_\nu\subseteq D((V_s)_{\neq 0}\otimes\allowbreak (V_s)_{\neq 0})$ for all $\nu\in \Phi_l$.

On the other hand, as $\Phi_s$ forms a single $W$-orbit of $\Phi$ and non-zero weight spaces in $V_l$ are one-dimensional, it suffices to prove that $D((V_s)_{\neq 0}\otimes\allowbreak (V_s)_{\neq 0})$ contains at least one weight vector of $V_l$ with a short root weight. There exist weight vectors $v_{\pm \mu_1}\in V_s$ of weight $\pm \mu_1\in \Phi_s$ with $D(v_{\mu_1}\otimes v_{-\mu_1}) = h \in \mathfrak{h}\setminus\{0\}$. By Lemma~\ref{lem:short-roots}, there is a short root $\mu_2\in \Phi_s$ with $\mu_2\neq \pm\mu_1$ such that $\mu_2(h)\neq 0$. As a consequence, we have $0\neq - \mu_2(h) e_{\mu_2} = [e_{\mu_2}, h] =  [e_{\mu_2},D(v_{\mu_1}\otimes v_{-\mu_1})] = D(\mathrm{ad}(e_{\mu_2}\otimes v_{\mu_1})\otimes v_{-\mu_1})+D(v_{\mu_1}\otimes \mathrm{ad}(e_{\mu_2}\otimes v_{-\mu_1}))$, and then it follows that $e_{\mu_2}\in D((V_s)_{\neq 0}\otimes\allowbreak (V_s)_{\neq 0})$. Hence, $(V_l)_\mu\subseteq D((V_s)_{\neq 0}\otimes\allowbreak (V_s)_{\neq 0})$ for all $\mu\in \Phi_s$, which completes the proof in the non-simply laced case.

The proof of $D((V_s)_{\neq 0}\otimes\allowbreak (V_s)_{\neq 0}) = V_l$ in the simply laced case follows from a similar argument.
\end{proof}

\begin{remark}
	For $\mathfrak{g}$ of type $B_n$ ($n\geq 2$), $D((V_s)_{\neq 0} \otimes (V_s)_{\neq 0})$ does not reach any weight vector whose weight is a short root, as the sum of two short roots can never be a short root. 
\end{remark}

\begin{theorem}\label{th:qm}
Let $\mathfrak{g}$ be at least of rank 2 and not of type $B_n\;(n\geq 2)$. Then any matrix coefficient of $V_s$ can be expressed as a polynomial in generalized minors of $V_s$ over $\mathbb{C}$.
\end{theorem}

\begin{proof}
If $\mathfrak{g}$ is simply laced, the proof of the theorem is straightforward and follows easily from Proposition \ref{prop:maps}(2). So, we take $\mathfrak{g}$ to be non-simply laced except for type $B_n$.

First, we prove that matrix coefficients of the form $C^{V_s}(z_i^*, v_\mu)$ can be written using generalized minors of $V_s$, where $z_i^* \in (V_s)^*_0$ is any functional of weight $0$ and $v_\mu \in V_s$ is a weight vector with weight $\mu\in \Phi_s$. By Proposition~\ref{prop:ad}(1), there exists an element $h\in (V_l)_0$ such that $v_\mu = \mathrm{ad}(h\otimes v_\mu)$. Therefore, we can write
\begin{equation}\label{eq1}
    C^{V_s}(z_i^*, v_\mu) = C^{V_s}(z_i^*, \mathrm{ad}(h \otimes v_\mu)) = C^{V_l\otimes V_s}(\mathrm{ad}^*(z_i^*), h \otimes v_\mu).
\end{equation}
By Proposition~\ref{prop:ad}(2), $\mathrm{ad}^*(z_i^*)$ must be a sum of terms of the form $f_{-\gamma} \otimes g_{\gamma}$, where $\gamma \in \Phi_s$, $f_{-\gamma} \in (V_l)^*_{-\gamma}$, and $g_{\gamma} \in (V_s)^*_{\gamma}$. So, we have
\begin{equation}\label{eq2}
    C^{V_s}(z_i^*, v_\mu)=\sum C^{V_l}(f_{-\gamma}, h) C^{V_s}(g_{\gamma}, v_\mu).
\end{equation}
Since $\gamma, \mu \in \Phi_s$, $C^{V_s}(g_{\gamma}, v_\mu)$ is by definition a generalized minor up to a scalar. Thus, it suffices to show that $C^{V_l}(f_{-\gamma}, h)$ can be written using generalized minors of $V_s$.

Since $V_l$ and $V_s$ are self-dual, we can identify $V_l$ and $V_s$ with $V_l^*$ and $V_s^*$ respectively. Let $\mathcal{D}\colon V_s^*\otimes V_s^* \to V_l^*$ be the map corresponding to $D$ under this identification.  Then, using Proposition~\ref{prop:maps}(2), we can write $f_{-\gamma}$ as $\mathcal{D}(\sum_k x^*_{\alpha_k} \otimes y^*_{\beta_k})$, where $x^*_{\alpha_k}$ and $y^*_{\beta_k}$ are functionals of non-zero weights. Substituting this, we get
\begin{equation}\label{eq3}
    C^{V_l}(f_{-\gamma}, h) = C^{V_l}\left(\mathcal{D}\left(\sum_k x^*_{\alpha_k} \otimes y^*_{\beta_k}\right), h\right) = \sum_k C^{V_s \otimes V_s}\left(x^*_{\alpha_k} \otimes y^*_{\beta_k}, \mathcal{D}^*(h)\right).
\end{equation}
Note that $\mathcal{D}^*(h)=\sum_{\delta \in  \Phi_s}c_\delta (v_{\delta}\otimes v_{-\delta})$ because $\mathcal{D}^*$ can be identified with $D^*$ and Proposition \ref{prop:maps}(2). Therefore,  $C^{V_l}(f_{-\gamma}, h)$ can be written in terms of generalized minors of $V_s$, so $C^{V_s}(z_i^*, v_\mu)$ is also a polynomial in generalized minors of $V_s$.
\vskip 0.1cm

Finally, let $z_j \in (V_s)_0$ and let $f \in (V_s)^*$ be either $z_i^*$ or $v_\mu^*$. By Proposition~\ref{prop:maps}(1), we can write $z_j = m(Z_j)$ where $Z_j = \sum_{\nu \in \Phi_s} a_\nu (v_\nu \otimes v_{-\nu})$ for some $a_\nu \in \mathbb{C}$. Writing $m^*(f) = \sum_r p_r^* \otimes q_r^*$ for some weight functionals $p_r^*, q_r^* \in (V_s)^*$, we obtain
\begin{equation}\label{eq4}
    C^{V_s}(f, z_j) = C^{V_s \otimes V_s}(m^*(f), Z_j) = \sum_{\nu \in \Phi_s} \sum_r a_\nu C^{V_s}(p_r^*, v_\nu) C^{V_s}(q_r^*, v_{-\nu}).
\end{equation}
Each factor on the right-hand side of \eqref{eq4} is either a generalized minor of $V_s$ up to a scalar or of the form $C^{V_s}(z_k^*, v_{\pm \nu})$ ($\nu \in \Phi_s$). We have already proved that the latter is a polynomial in generalized minors of $V_s$, thereby completing the proof.
\end{proof}

As a consequence, we immediately obtain our main result.

\begin{cor}[= Theorem~\ref{th:main}]
   Let $G$ be any connected simply connected simple complex Lie group. The coordinate ring $\mathcal{O}(G)$ is generated by generalized minors over $\mathbb{C}$. 
\end{cor}

\begin{proof}
     For $G$ of the types $A_n, B_n, C_n, D_n, E_6$, and $E_7$, $G$ admits a tensor generating family of minuscule representations, so $\mathcal{O}(G)$ is generated by their matrix coefficients. Since every weight of a minuscule module is extremal, all of its weight-basis matrix coefficients are generalized minors up to a scalar, thus proving the desired inclusion. For $G$ of the types $G_2, F_4$, and $E_8$, $\mathcal{O}(G)$ is generated by matrix coefficients of $V_s$, and Theorem \ref{th:qm} then completes the proof.
\end{proof}

\begin{remark}
Write $A:=\mathbb{Z}[q,q^{-1}]$ and assume that $A$ acts on $\mathbb{Z}$ by sending $q \mapsto 1$. Let $\mathcal{O}^A_q(G)$ be Lusztig's $A$-form of the quantized coordinate ring $\mathcal{O}_q(G)$ over $\mathbb{Q}(q)$. When $G$ is not of type $F_4$, Oya et al. in \cite{OQY} recently proved a sharp integral version of the above result. That is, the $\mathbb{Z}$-form $\mathcal{O}_{\mathbb{Z}}(G):= \mathbb{Z}\otimes_{A}\mathcal{O}^{A}_q(G)$ is generated by generalized minors over $\mathbb{Z}$ when $G$ is not of type $F_4$ or $G_2$, and over $\mathbb{Z}[\frac{1}{2}]$ when $G$ is of type $G_2$. But such an integral refinement is not yet known to us in type $F_4$.
\end{remark}

\begin{cor}
     The coordinate ring $\mathcal{O}(G)$ coincides with the Berenstein--Fomin--Zelevinsky cluster algebra and its upper cluster algebra, or in other words, $\mathcal{O}(G)=\mathcal{A}=\mathcal{U}$. 
\end{cor}

\begin{proof}
By \cite[Theorem 1.4]{GLS13} and \cite[Theorem 3.1]{Oya}, proving the equalities $\mathcal{O}(G)=\mathcal{A}=\mathcal{U}$ reduces to showing $\mathcal{O}(G)\subseteq \mathcal{A}$, which is then supplied by Theorem~\ref{th:main}, completing the proof.
\end{proof}

\subsection*{Acknowledgements} We would like to thank Hironori Oya for his feedback and comments. The second author is supported by the Indian Statistical Institute through a PhD fellowship.

\subsection*{AI declaration} LLMs were used for literature search and during the writing phase, as an aid to find errors, and typesetting some parts of the LaTeX manuscript. The underlying ideas of the main results, including the proof strategy, are due to the authors.

\end{document}